\documentstyle[12pt]{article}

\def\beq{\begin{equation}}
\def\ee{\end{equation}}
\def\bib#1{[{\ref{#1}}]}
\def\psiba{\overline{\psi}}
\def\vphiba{\overline{\varphi}}
\def\zba{\overline{z}}
\def\paz{\partial_{z}}
\def\CC{\mbox{\boldmath $C$}}
\def\RR{R}
\def\SS{S}
\def\CP{CP}
\def\GG{G}
\def\HH{H}
\def\QQ{Q}
\def\pazba{\partial_{\overline{z}}}
\newtheorem{theorem}{Theorem}[section]
\newtheorem{proposition}{Proposition}
\newtheorem{corollary}{Corollary}
\begin{document}           


            \title{
Generalized Weierstrass representation for surfaces in
multidimensional Riemann spaces}
 \author{B.G. Konopelchenko and G. Landolfi \\
{\em Dipartimento di Fisica, Universit\'{a} di Lecce, 
73100 Lecce, Italy}\\ and \\
{\em I.N.F.N., Sezione di Lecce, 73100 Lecce, Italy} } 

 \maketitle

      \begin{abstract}

Generalizations of the Weierstrass formulae to generic surface
immersed into $\RR^4$, 
$\SS^4$ and into multidimensional Riemann spaces are proposed. 
Integrable deformations of surfaces in these spaces via the 
modified Veselov-Novikov equation are discussed.
      \end{abstract}


\section{Introduction}
\setcounter{equation}{0}

Theory of immersion and deformations of surfaces has been 
important 
part of the classical differential geometry (see {\em e.g.} 
\bib{r1}-\bib{r3}). 
Various methods to describe immersions and different types of 
deformations 
have been considered. New results in this field have 
been discussed, for instance, in \bib{r4}-\bib{r5}). 

Surfaces and their dynamics are key 
ingredients in a number of phenomena in 
physics too (see {\em e.g.} \bib{r6}-\bib{r8}). 
They are, for instance, surface waves, 
propagation of flame fronts, growth of crystals, deformation of 
membranes, 
dynamics of vortex sheets, many problems of hydrodynamics 
connected with motion of boundaries between 
region of differing densities and viscosities. 
Number of papers has been devoted to a study and application 
of the integrals over surfaces in gauge field theories, string 
theory, quantum gravity and statistical physics 
\bib{r6}-\bib{r8}).

Direct approaches to describe surfaces always have been of great 
interest. 
The classical Weierstrass formulae for minimal surfaces immersed
in the three-dimensional Euclidean space 
$\RR^3$ is the best known example of such 
an approach. Only recently \bib{r9}-\bib{r10}-\bib{r11}
 the Weierstrass formulae have been 
generalized to the case of generic surfaces in $\RR^3$. 
During the last two 
years the generalized Weierstrass formulae have been used 
intensively to study both global properties of surfaces 
in $\RR^3$ and their integrable deformations. 

In this paper we present the generalizations of the Weierstrass 
representation for surfaces immersed into the multi-dimensional
Euclidean 
and Riemann spaces. The cases of the four-dimensional Euclidean
space $\RR^4$ and space $\SS^4$ 
of constant curvature are considered in  detail. The 
comparison of our Weierstrass formulae for surfaces in 
$\RR^4$ and those of 
conformal immersion in $\RR^4$ is given. 
The properties of the Willmore 
functional for the immersion in $\RR^4$ and $\SS^4$ 
are studied. The 
Weierstrass representations for immersion into pseudo-euclidean 
spaces with signatures $(+,+,+,-)$ and $(+,+,-,-)$ are presented.
Surfaces on Lie groups are discussed too. 

Integrable deformations of surfaces via the modified 
Veselov-Novikov 
equation are considered. It is shown that the Willmore 
functional (or Helfrich-Polyakov action) 
is invariant under such deformations.

The paper is organized as follows. In section $2$ a brief review 
of the old and generalized Weierstrass formulae for surfaces in 
$\RR^3$ 
is given. The Weierstrass representation for generic surfaces 
immersed into $\RR^4$ is derived in section $3$. 
An extension to four-dimensional Riemann spaces, in particular 
to $\SS^4$ and the Minkowski space, is presented in section 
$4$. The Weierstrass representations 
for surfaces in multidimensional 
spaces are discussed in section $5$. 
Integrable deformations are considered in section $6$.

\section{The old and generalized Weierstrass formulae for 
surfaces in $\RR^3$}
\setcounter{equation}{0}

Here for the sake of convenience we will briefly remind the 
old Weierstrass representation for minimal surfaces in $\RR^3$ 
and will review recent results 
concerned its generalization to generic surfaces in $\RR^3$. 

An original Weierstrass formulae (see {\em e.g.} 
\bib{r1}-\bib{r3}) starts with two arbitrary holomorphic 
functions $\varphi(z)$ and $\chi(z)$ of the complex
variables $z$, $\zba \in \CC$. Then one 
introduces the three quantities $X^i(z,\zba)$ ($i=1,2,3$) 
as follows:
\begin{eqnarray}
X^1 & = & Re 
\left[ i \int_{\Gamma} {\left( \varphi^2 +\chi^2 \right)dz'}
\right] \; ,  \nonumber \\
X^2 & = & Re \left[ 
\int_{\Gamma}{\left( \varphi^2 -\chi^2 \right)dz'}
\right] \; , \nonumber \\
X^3 &  = & - Re 
\left[ 2 \int_{\Gamma} { \varphi \chi dz'}\right] \;\; .
\label{eq:2.1}
\end{eqnarray}
Finally one treats $X^i(z,\zba)$ ($i=1,2,3$) as the 
coordinates of a 
surface in $\RR^3$: this surface is a {\em minimal} 
one ({\em i.e.} has zero 
mean curvature) and the parametric lines $z=const.$ and 
$\zba=const$ 
are the minimal lines. The Weierstrass formulae have been the 
most powerful tool in study of minimal surfaces. 

An analog of the Weierstrass formulae for surfaces of 
prescribed (non zero) 
mean curvature have been proposed by Kenmotsu in $1979$ 
\bib{r9}. The Kenmotsu representation is given by
\begin{eqnarray}
X^i &  = & Re \left[ \int^{z}{ \eta \phi^i dz'} \right]
\label{eq:2.2}
\end{eqnarray}
where $\vec{\phi}=[1-f^2, i(1+f^2), 2f]$ 
and the functions $f$ and $\eta$ 
obey the following compatibility condition
\beq
\left( \log \eta \right) _{\zba} =
- \frac{ 2 \overline{f} f_{\zba} }{ 1+|f|^2 } \;\; .
\label{eq:2.3}
\ee
Here and below the bar denotes the complex conjugation. Then 
the mean curvature $H$ is
\beq
H=-\frac{2\overline{f}_z}{\eta \left( 1+|f|^2 \right)^2} \;\; .
\label{eq:2.4}
\ee
It was proved in \bib{r9} that any surface in $\RR^3$ 
can be represented in 
such a form. This representation of surfaces deals 
basically with the 
Gauss map for generic surface in $\RR^3$ 
(see also \bib{r12}-\bib{r13}). 

Another generalization of the Weierstrass formulae to generic 
surfaces in $\RR^3$ 
have been proposed independently by one of the authors in 
$1993$ (see \bib{r10} and \bib{r11}). 
It starts with the linear system (two-dimensional Dirac equation)
\beq
\begin{array}{l}
\psi_{z}=p \varphi \;\; ,\\
\varphi_{\zba}=-p \psi
\end{array}
\label{eq:2.5}
\ee
where $\psi$ and $\varphi$ are complex-valued functions of 
$z$,$\zba \in \CC$ and $p(z,\zba)$ is a real-valued function. 
Then one defines the three real-valued functions $X^1(z,\zba)$,
$X^2(z,\zba)$ 
and $X^3(z,\zba)$ by the formulae \bib{r10}-\bib{r11}
\begin{eqnarray}
X^1 & + & i X^2  
= i \int_{\Gamma} {\left( \overline{\psi}^2 dz'-
\overline{\varphi}^2 d\zba' \right)} \; ,  \nonumber \\
X^1 & - &  i X^2  = i \int_{\Gamma}{\left( \varphi^2 dz'-
\psi^2 d\zba' \right)}  \; ,\nonumber \\
X^3 &  = & - \int_{\Gamma} { \left( \overline{\psi} \varphi dz'
+\psi \overline{\varphi} d\zba' \right)} \;\; .
\label{eq:2.6}
\end{eqnarray}
where $\Gamma$ is an arbitrary curve in $\CC$. 
In virtue of (\ref{eq:2.5}) 
the {\em r.h.s.} in (\ref{eq:2.6}) do not depend on the 
choice of $\Gamma$. 
If one now treats $X^i(z,\zba)$ as the coordinates in 
$\RR^3$ then the 
formulae (\ref{eq:2.5}), (\ref{eq:2.6}) define a conformal 
immersion of surface into $\RR^3$ with the induced metric of 
the form
\beq
ds^2=u^2 dz d\zba \;\; , \;\; u=|\psi|^2+|\varphi|^2 \;\; ,
\label{eq:2.7}
\ee
with the Gauss curvature
\beq
K=-\frac{4}{u^2} \left[ \log {u} \right]_{z\zba}
\label{eq:2.8}
\ee
and the mean curvature
\beq
H=2\frac{p}{u} \; \; .
\label{eq:2.9}
\ee
At $p=0$ one gets minimal surfaces and the formulae 
(\ref{eq:2.6}) are reduced to 
the old Weierstrass formulae (\ref{eq:2.1}) 
under the identification $\varphi=\varphi$, $\chi=\psiba$. 

It turned out that the Kenmotsu formulae (\ref{eq:2.2}),
(\ref{eq:2.3})
and the generalized Weierstrass formulae 
(\ref{eq:2.5}),(\ref{eq:2.6}) are 
equivalent to each other. The relation between the functions 
$(f, \eta)$ and $(\psi, \varphi)$ is the following \bib{r14}
\begin{eqnarray}
f=i\frac{\psiba}{\varphi}\;\; , \;\; \eta=i \varphi^2
\label{eq:2.10}
\end{eqnarray}
and
\beq
p= -\frac{ \eta f_{\zba} }{ \sqrt{\eta 
\overline{\eta} } \left( 1+|f|^2 \right) } \; \; .
\label{eq:2.11}
\ee
So, all results proved for the Kenmotsu formulae \bib{r9} 
and associated Gauss map \bib{r13} 
in $\RR^3$ are valid also for the generalized Weierstrass 
formulae (\ref{eq:2.6}). In particular, it implies 
immediately that any surface in $\RR^3$ 
can be represented via (\ref{eq:2.5}), (\ref{eq:2.6}).

Though the representations (\ref{eq:2.2})-(\ref{eq:2.4}) 
and (\ref{eq:2.5}), (\ref{eq:2.6}) are equivalent, the latter 
provides us certain advantages. 
They are mainly due to the fact that in the generalized 
Weierstrass 
formulae the functions $\psi$ and $\varphi$ obey linear 
equations (\ref{eq:2.5}) while for the Kenmotsu 
formulae the nonlinear constraint (\ref{eq:2.3}) 
is hard to treat. This circumstance had allowed to simplify 
essentially an analysis that had lead to several interesting 
results both of local and global character \bib{r14}-\bib{r21}. 
It occurred, in particular, that the Willmore functional 
(see {\em e.g.} \bib{r22}) or the Helfrich-Polyakov 
action (see \bib{r6}-\bib{r8}) $W=\int {\vec{H}^2 \left[ 
dS \right] }$ has a 
very simple form: $W=4 \int{ p^2 dx \; dy}$  ($z=x+iy$)
\bib{r14}-\bib{r15}. 

One of the advantages of the generalized Weierstrass formulae 
(\ref{eq:2.5}), (\ref{eq:2.6}) 
is that they allow to construct a new class of deformations 
of surfaces via the modified Veselov-Novikov equation 
\bib{r10}-\bib{r11}. The 
characteristic feature of these integrable deformations is 
that the Willmore functional 
remains invariant \bib{r14}-\bib{r15}. Thus, the generalized 
Weierstrass representation (\ref{eq:2.6}) has been 
proved to be an effective tool to study 
surfaces in $\RR^3$ and their integrable deformations.

We would like to emphasize that the idea to 
generate surfaces via solutions 
of linear equations is, in fact, the old idea of the classical 
differential 
geometry (see discussion in \bib{r11}). In \bib{r3} 
one can find the two representations 
of these type in addition to the Weierstrass formulae. The 
first is given by the Lelieuvre's formula which is well 
known in affine geometry (see {\em e.g.} \bib{r23}). 
Another example 
(\bib{r3}, p. $82$) is provided by the equation
\beq
\theta_{\xi\eta}-\left( \log \lambda \right)_\eta \theta_\xi
-\lambda^2 \theta=0
\label{eq:2.13}
\ee
where $\xi$, $\eta$ are real variables and $\lambda$ is 
a real-valued function. 
It is stated in \bib{r3} that two solutions of 
(\ref{eq:2.13}) define, via certain integral formulae,
 a surface in $\RR^3$ parameterized by 
minimal lines, but no calculation of the metric and 
curvature is given. 
This example, seems, was forgotten completely until it 
had been found during the preparation of the second 
paper \bib{r11} on the generalized Weierstrass formulae. 
The representation (\ref{eq:2.13}) 
is rather close to that of (\ref{eq:2.5}), 
(\ref{eq:2.6}). Indeed, equation 
(\ref{eq:2.13}) can be rewritten as the system
\beq
\begin{array}{l}
\theta_{\xi}=\lambda \varphi \; , \\
\varphi_{\eta}=\lambda \theta
\end{array}
\label{eq:2.14}
\ee
where $\varphi$ is a new function. If one takes two solutions 
$(\theta, \varphi)$ 
and $(\tilde{\theta}, \tilde{\varphi})$ 
of the system (\ref{eq:2.14}) then 
the formulae given in \bib{r3} (p. $82$) take the form
\begin{eqnarray}
X^1 & + & i X^2  = \int{\left( \theta^2 d\eta +\varphi^2 d\xi
\right)} \; , \nonumber \\
X^1 & - & i X^2  = \int{\left( \tilde{\theta}^2 d\eta +
\tilde{\varphi}^2 d\xi \right)} \; , \nonumber \\
X^3 & = & i\int{\left( \theta \tilde{\theta} d\eta +
\varphi \tilde{\varphi} d\xi \right)} \;\; . 
\label{eq:2.15}
\end{eqnarray}
However, in contrast to the representation 
(\ref{eq:2.5})(\ref{eq:2.6}), 
the formulae (\ref{eq:2.14}), 
(\ref{eq:2.15}) do not define a real 
surface in $\RR^3$.

We would like to note that some results in \bib{r24} and 
\bib{r25} were close to 
the generalized Weierstrass representation (\ref{eq:2.6}). 
In \bib{r24} a formula similar 
to (\ref{eq:2.5}) for constant mean curvature surfaces has been 
discussed. In \bib{r25} the system (\ref{eq:2.5}) had appeared 
within the quaternionic description of surfaces in $\RR^3$ 
(formula ($2.19$) of \bib{r25}). 
But in \bib{r25} it was accompanied by another 
two equations (equation ($2.16$) of \bib{r25}) 
which are indispensable in the Sym's type approach. So the 
meaning of the system (\ref{eq:2.5}), seems, had been missing.
The generalized Weierstrass type formulae admit also a beautiful 
formulation within the spinor representations of surfaces 
\bib{r26}-\bib{r27}.
   
\section{The Weierstrass representation for immersion into 
$\RR^4$}
\setcounter{equation}{0}

An extension of the representation (\ref{eq:2.5}), (\ref{eq:2.6}) 
to the four dimensional Euclidean space is as follows. 
Let $\psi_1$, $\varphi_1$ 
and $\psi_2$, $\varphi_2$ be two independent solutions 
of the system 
(\ref{eq:2.5}), {\em i.e.}
\begin{eqnarray}
\begin{array}{l}
\psi_{1 z}=p \varphi_1 \;\; ,\\
\varphi_{1 \zba}=-p \psi_1
\end{array}
&   & 
\begin{array}{l}
\psi_{2 z}=p \varphi_2 \; , \\
\varphi_{2 \zba}=-p \psi_2 \;\;\; .
\end{array}
\label{eq:3.1}
\end{eqnarray}
Equations (\ref{eq:3.1}) imply that
\begin{eqnarray}
\left( \psi_1 \psi_2 \right)_{z} = -\left( \varphi_1 \varphi_2 
\right)_{\zba}
\;\; , \;\; 
\left( \psi_1 \vphiba_2 \right)_{z}=\left( \varphi_1 \psiba_2 
\right)_{\zba} \;\; .
\label{eq:3.2}
\end{eqnarray}
As a consequence there are four functions 
$X^i(z,\zba)$ ($i=1,2,3,4$) such that
\begin{eqnarray}
dX^1 & = & \frac{i}{2} \left( \psiba_1 \psiba_2 +
\varphi_1 \varphi_2 \right) dz +c.c \; , \nonumber \\
dX^2 & = & \frac{1}{2} \left( \psiba_1 \psiba_2 -
\varphi_1 \varphi_2 \right) dz +c.c \; , \nonumber \\
dX^3 & = & -\frac{1}{2} \left( \psiba_1 \varphi_2 +
\psiba_2 \varphi_1 \right) dz + c.c. \; , \nonumber \\
dX^4 & = & \frac{i}{2} \left( \psiba_1 \varphi_2 -
\psiba_2 \varphi_1 \right) dz + c.c. \;\; .
\label{eq:3.3}
\end{eqnarray}
where $c.c.$ means complex conjugated previous term.
We treat now these functions $X^i(z,\zba)$ as the coordinates 
of surface in $\RR^4$. For components of induced metric
\begin{eqnarray}
g_{zz}=\sum_{i=1}^4 {\left( X^i_z \right)^2}
=\overline{g_{\zba \; \zba}} 
& , & 
g_{z\zba}=\sum_{i=1}^4 { \left( X^i_z X^i_{\zba} \right)}
\label{eq:3.4}
\end{eqnarray}
one gets
\beq
g_{zz}=g_{\zba \; \zba}=0
\label{eq:3.5}
\ee
and
\beq
g_{z \zba}=\frac{1}{2} \left( |\psi_1|^2+|\varphi_1|^2 \right)
\left( |\psi_2|^2+|\varphi_2|^2 \right)
\;\; .
\label{eq:3.6}
\ee
Further, two normal vectors $\vec{N}_1$, $\vec{N}_2$ are
\begin{equation}
\vec{N}_1=\sqrt{\frac{|\varphi_1|^2 |\varphi_2|^2}{u_1 u_2}} 
Re \left(\vec{A} \right) \;\;, \;\;
\vec{N}_2=\sqrt{\frac{|\varphi_1|^2 |\varphi_2|^2}{u_1 u_2}}
Im \left(\vec{A} \right) \;\; .
\label{eq:3.7}
\ee
where 
\begin{eqnarray}
u_k & = & \left( |\psi_k|^2+|\varphi_k|^2 \right) \;\;\;\; , 
\;\;\;\; k=1,2 \\
\vec{A} & = & 
\left[ 
i \left( \frac{\psi_1}{\vphiba_1} - \frac{\psiba_2}{\varphi_2}
\right), -\frac{\psi_1}{\vphiba_1}-\frac{\psiba_2}{\varphi_2}, 
1-\frac{\psi_1 \psiba_2}{\vphiba_1 \varphi_2}, 
-i\left(  1+\frac{\psi_1 \psiba_2}{\vphiba_1 \varphi_2}
\right) \right] \, .
\label{eq:3.8}
\end{eqnarray}
The mean curvature vector 
$\vec{H}=\frac{\vec{X}_{z\zba}}{g_{z\zba}}$ 
is given by
\begin{eqnarray}
\vec{H} & =  & \frac{2p}{u_1 u_2} 
Re\left[ -i \left(\psi_1 \varphi_2 +\psi_2 \varphi_1 \right), 
\left(\psi_1 \varphi_2 +\psi_2 \varphi_1 \right) , \right. 
\nonumber \\
& & \left. \;\;\;\;\;\;\;\;\;\;\; \;\;
\left( \psi_1 \psiba_2 -\varphi_1 \vphiba_2 \right), 
i \left( \psi_1 \psiba_2 -\varphi_1 \vphiba_2 \right)
 \right] \; .
\label{eq:3.9}
\end{eqnarray}
The components $h_1$, $h_2$ of $\vec{H}$ along $\vec{N}_1$ 
and $ \vec{N}_2 $
( $\vec{H}=h_1\vec{N}_1+h_2\vec{N}_2 $ ) are
\begin{eqnarray}
h_1 & = & -p  
\frac{ \varphi_1 \vphiba_2 + 
\vphiba_1 \varphi_2}{ \sqrt{ u_1 u_2 
|\varphi_1|^2 |\varphi_2|^2}} \; , \nonumber \\
h_2 & = & ip  
\frac{ \varphi_1 \vphiba_2 - 
\vphiba_1 \varphi_2}{ \sqrt{ u_1 u_2 
|\varphi_1|^2 |\varphi_2|^2}} \;\; .
\label{eq:3.10} 
\end{eqnarray}
So, the mean curvature $H=\sqrt{ \sum_{i=1}^{4} {H^i H^i} }=
\sqrt{h^2_1+h^2_2}$ is equal to
\beq
H=\frac{2p}{\sqrt{u_1 u_2}} \;\; .
\label{eq:3.12}
\ee
Then the Gaussian curvature is
\beq
K=-\frac{2}{ u_1 u_2 } \left[ \log {\left( u_1 u_2 \right)} 
\right]_{z\zba} \;\;. 
\label{eq:3.13}
\ee
Finally, the Willmore functional $W=\int {\vec{H}^2 \left[ dS
\right] }$ is given by 
\beq
W=4 \int {p^2 dx dy} \;\; .
\label{eq:3.14}
\ee
Thus, we have the following
\begin{theorem}
 The generalized Weierstrass formulae
\begin{eqnarray}
X^1 & = & \frac{i}{2} \int_{\Gamma} {\left[ \left( 
\psiba_1 \psiba_2+\varphi_1 \varphi_2 \right) dz'-
\left( \psi_1 \psi_2 +\vphiba_1 \vphiba_2 \right) d\zba' 
\right]}  \; , \nonumber \\
X^2 & = & \frac{1}{2} \int_{\Gamma} {\left[ \left( 
\psiba_1 \psiba_2-\varphi_1 \varphi_2 \right) dz' +
\left( \psi_1 \psi_2 - \vphiba_1 \vphiba_2 \right) d\zba' 
\right]}  \; , \nonumber \\
X^3 & = & -\frac{1}{2} \int_{\Gamma} {\left[ \left( 
\psiba_1 \varphi_2+\psiba_2 \varphi_1 \right) dz'+
\left( \psi_1 \vphiba_2 +\psi_2 \vphiba_1  \right) d\zba' 
\right]}  \; , \nonumber \\
X^4 & = & \frac{i}{2} \int_{\Gamma} {\left[ \left( 
\psiba_1 \varphi_2 - \psiba_2 \varphi_1 \right) dz'-
\left( \psi_1 \vphiba_2 -\psi_2 \vphiba_1  \right) d\zba' 
\right]}  
\label{eq:3.15}
\end{eqnarray}
where 
\begin{eqnarray}
\begin{array}{l}
\psi_{\alpha z}=p \varphi_{\alpha} \; , \\
\varphi_{\alpha \zba}=-p \psi_\alpha
\end{array}
 &, & \alpha=1,2
\label{eq:3.16}
\end{eqnarray}
$p(z,\zba)$ is a real valued function, $\Gamma$ is a contour in 
$\CC$, define the conformal immersion of a surface into $\RR^4$. 
The induced metric is of the form
\beq
ds^2=u_1 u_2 dz d\zba
\label{eq:3.17}
\ee
where $u_\alpha=|\psi_\alpha|^2+
|\varphi_\alpha|^2$ ($\alpha=1,2$), the Gaussian 
and mean curvatures are
\begin{eqnarray}
K=-\frac{2}{u_1 u_2} 
\left[ \log { \left( u_1 u_2 \right) } \right] _{z \zba} 
\;\; , \;\; H=\frac{2p}{\sqrt{u_1 u_2}} \;\; .
\label{eq:3.18}
\end{eqnarray}
The total squared mean curvature (Willmore functional) 
is given by
\beq
W=4\int {p^2 dx dy} \;\; .
\label{eq:3.20}
\ee
\end{theorem}
The generalized Weierstrass 
representation (\ref{eq:3.15}) defines surface 
in $\RR^4$ up to translations. In the particular case 
$\psi_2=\pm \psi_1$, $\varphi_2=\pm\varphi_1$, 
$X^4_z=X^4_{\zba}=0$ and the formulae 
(\ref{eq:3.16})-(\ref{eq:3.20}) are reduced to those 
(\ref{eq:2.5})-(\ref{eq:2.9}) of $\RR^3$ case. 
\begin{corollary}
 Minimal surfaces in $\RR^4$ are given by the Weierstrass 
representation (\ref{eq:3.15}), (\ref{eq:3.16}) 
with $p(z,\zba)=0$. 
Surfaces of constant mean curvature $H$ are given by the 
formulae (\ref{eq:3.15}) where $\psi_\alpha$, $\varphi_\alpha$ 
($\alpha=1,2$) obey the system of equations
\begin{eqnarray}
\psi_{\alpha z} & = & \frac{H}{2} 
\sqrt{\left(|\psi_1|^2+|\varphi_1|^2 \right)
\left(|\psi_2|^2+|\varphi_2|^2 \right)} \; \varphi_\alpha 
\;\; , \nonumber \\
\varphi_{\alpha \zba} & = & -\frac{H}{2}
\sqrt{\left(|\psi_1|^2+|\varphi_1|^2 \right)
\left(|\psi_2|^2+|\varphi_2|^2 \right)} \; \psi_\alpha 
\;\; .
\label{eq:3.21}
\end{eqnarray}
\end{corollary}
At $\psi_2=\pm \psi_1$, $\varphi_2=\pm\varphi_1$ the system 
(\ref{eq:3.21}) is reduced to a simpler one \bib{r11} which has 
been studied in \bib{r16}.

Note the equations (\ref{eq:3.3}) can be represented in the form
\begin{eqnarray}
d\left( X^1+iX^2 \right) & = & i\psiba_1 \psiba_2 dz -  
i \vphiba_1 \vphiba_2 d\zba \; ,\nonumber \\
d\left( X^1-iX^2 \right) & = & i \varphi_1 \varphi_2 dz  
- i\psi_1 \psi_2 d\zba \; , \nonumber \\
d\left( X^4+iX^3 \right) & = & -i\psiba_2 \varphi_1 dz -
i \psi_1 \vphiba_2 d\zba \; , \nonumber \\
d\left( X^4-iX^3 \right) & = & i\psiba_1 \varphi_2 dz +
i \psi_2 \vphiba_1 d\zba
\label{eq:3.22}
\end{eqnarray}
which reveals a symmetry between the pairs of coordinates 
$(X^1,X^2)$, $(X^3,X^4)$. 

The formulae (\ref{eq:3.3}) can be also rewritten in a spinor 
representation type form
\beq
d\left( \sigma_1 X^2+\sigma_2X^1 - \sigma_3 X^3 +
i I X^4 \right) = V_2^{\dagger} 
\left( 
\begin{array}{cc}
0 & dz \\
d\zba & 0
\end{array}
\right)
V_1
\label{eq:3.23}
\ee
where 
\beq
V_{1,2}= \left(
\begin{array}{cc}
\psi_{(1,2)} & -\vphiba_{(1,2)} \\
\varphi_{(1,2)} & \psiba_{(1,2)}
\end{array}
\right) \; ,
\nonumber
\ee
$\sigma_i$ ($i=1,2,3$) are the Pauli matrices and $I$ is the 
identity matrix. 

The condition (\ref{eq:3.5}), that an immersion is conformal, 
written as
\beq
\left( X^1_z\right)^2+\left( X^2_z\right)^2+
\left( X^3_z\right)^2+\left( X^4_z\right)^2=0
\label{eq:3.24}
\ee
defines the complex quadric $Q_2$
\beq
w_1^2+w_2^2+w_3^2+w_4^2=0
\label{eq:3.25}
\ee
in $\CP^3$ where $w_i$ ($i=1,2,3,4$) are homogeneous coordinates. 
A diffeomorphism of ${\QQ}_2$ to the Grassmannian 
$\GG_{2,4}$ of oriented 
$2-$planes in $\RR^4$ allows us to define the Gauss map $G(z)$ 
for a surface represented by the 
generalized Weierstrass formulae (\ref{eq:3.15}). It is given by
\beq
\vec{G}(z)=
\left[ i \left( \psiba_1 \psiba_2 +\varphi_1 \varphi_2 \right),
\psiba_1 \psiba_2 -\varphi_1 \varphi_2, 
-\psiba_1 \varphi_2 - \psiba_2 \varphi_1, 
i \left( \psiba_1 \varphi_2 - \psiba_2 \varphi_1 \right)
\right] \;\; .
\label{eq:3.26}
\ee
The Gauss map for surfaces immersed 
into $\RR^4$ has been studied 
earlier in the paper \bib{r13}. In \bib{r13} the Gauss map 
$\vec{G}(z)$ has been parameterized as follows
\beq
\vec{G}(z)= 
\left[ \left(1+f_1f_2 \right), i\left(1-f_1f_2 \right), 
\left(f_1-f_2 \right), -i \left(f_1 + f_2 \right) \right]
\label{eq:3.27}
\ee
where $f_1$ and $f_2$ are complex-valued functions. A surface 
in $\RR^4$ is then defined by \bib{r13}
\beq
\vec{X}=\int^{z} {Re\left(\eta \vec{G} dz \right) }
\label{eq:3.28}
\ee
where $f_1$ and $f_2$ satisfy the compatibility conditions
\beq
Im \left[ 
\left( \frac{f_{1 z\zba} }{f_{1\zba} } -
\frac{2\overline{f}_1 f_{1z} }{1+|f_1|^2} \right)_{\zba} +
\left( \frac{f_{2 z\zba}}{f_{2\zba}}-
\frac{2\overline{f}_2 f_{2z} }{1+|f_2|^2} \right)_{\zba}
\right] = 0
\label{eq:3.29}
\ee
and
\beq
|F_1|=|F_2|
\label{eq:3.30}
\ee
where $ F_{i} =  f_{i \zba} \left( 1+|f_i|^2 \right)^{-1} $, 
$i=1,2$. The function $\eta$ is given by 
\beq
\overline{\eta}^2=-\frac{4 F_1 F_2}{H^2 \left( 1+|f_1|^2 \right) 
\left( 1+|f_2|^2 \right) }
\label{eq:3.31}
\ee
where the mean curvature $H$ is expressed via $f_1$ and $f_2$ by
\beq
2 \left( \log H \right)_z = 
 \frac{f_{1 z\zba}}{f_{1\zba}}-
\frac{2\overline{f}_1 f_{1z} }{1+|f_1|^2}
+\frac{f_{2 z\zba}}{f_{2\zba}}-
\frac{2\overline{f}_2 f_{2z} }{1+|f_2|^2} \;\; .
\label{eq:3.32}
\ee
Similar to the three-dimensional case this representation 
includes the complicated compatibility conditions.
\begin{theorem}
 The generalized Weierstrass representation 
(\ref{eq:3.15})-(\ref{eq:3.20}) implies the 
Gauss map type representation 
(\ref{eq:3.28})-(\ref{eq:3.32}) via the substitution
\begin{eqnarray}
\eta=i \varphi_1 \varphi_2
\;\; , \;\;
f_1=i \frac{ \psiba_1 }{ \varphi_1 }
\;\; , \;\;
f_2=-i\frac{ \psiba_2 }{ \varphi_2 } \;\; . 
\label{eq:3.33}
\end{eqnarray}

\end{theorem}

The proof is straightforward: equations (\ref{eq:3.16})
and (\ref{eq:3.33}) give the constraints 
(\ref{eq:3.29})-(\ref{eq:3.31}) with
\beq
p=-iF_1 \frac{\varphi_1}{\vphiba_1}
=iF_2 \frac{\varphi_2}{\vphiba_2}
\label{eq:3.34}
\ee
while (\ref{eq:3.15}) is converted into (\ref{eq:3.28}).

\section{Surfaces in four-dimensional Riemann space}
\setcounter{equation}{0}

The results of the previous section apparently can be 
extended to the case 
of immersion into generic four-dimensional Riemann space 
with the metric tensor $g_{ik}$. 
\begin{proposition}
 The formulae (\ref{eq:3.15}), (\ref{eq:3.16}) define an 
immersion of surface into the four-dimensional Riemann space 
with the metric tensor 
$g_{ik}$. The induced metric is
\beq
ds^2=g_{zz}dz^2+2g_{z\zba}dz d\zba +g_{\zba \zba} d\zba^2
\label{eq:4.1}
\ee
where
\begin{eqnarray}
g_{zz}=g_{ik} X^i_z X^k_z =\overline{g_{\zba \; \zba} }\;\; , 
\;\;
g_{z\zba}=g_{ik} X^i_z X^k_{\zba} \;\; .
\label{eq:4.2}
\end{eqnarray}
\end{proposition}
The Gaussian and mean curvature are calculated 
straightforwardly.

In the case of conformally-Euclidean 
spaces $g_{ik}=e^{2\sigma}\delta_{ik}$ 
($i,k=1,2,3,4$, $\sigma$ is a function and $\delta_{ik}$ is the 
Kronecker symbol) the immersion is the conformal one:
\beq
ds^2=e^{2\sigma} u_1 u_2 dz d\zba \;\; .
\label{eq:4.3}
\ee
The Gaussian and mean curvatures are 
\beq
K=-2 e^{-2\sigma} \frac{ \left[ 2\sigma +\log {\left( u_1u_2 
\right)} \right]_{z\zba}}{u_1 u_2} \;\; , \;\;
H=2 e^{-\sigma} \frac{p}{\sqrt{u_1 u_2}} \;\; .
\label{eq:4.4-4.5}
\ee
For the Willmore functional one gets
\beq
W=4\int {p^2 dx dy} \;\; .
\label{eq:4.6}
\ee
A special case of immersions into the space $\SS^4$ of constant 
curvature had attracted recently the particular interest (see 
{\em e.g.} \bib{r22},\bib{r28}). To 
describe it we choose the Riemann form for the metric of 
$\SS^4$, {\em i.e.} (see {\em e.g.} \bib{r29})
\beq
e^{\sigma}=\left[ 1+\frac{K_0}{4} 
\sum_{i=1}^4\left(X^i \right)^2 \right]^{-1}
\label{eq:4.7}
\ee
where $K_0$ is the curvature. Then the formulae 
(\ref{eq:3.15}), (\ref{eq:3.16}), (\ref{eq:4.3})-(\ref{eq:4.6})
define the conformal immersion of a surface into $\SS^4$. At
$\psi_2=\pm\psi_1$, $\varphi_2=\pm\varphi_1$ ($X^4=0$) 
one has the conformal immersion into $\SS^3$. 
The generalized Weierstrass representation provides 
us an effective method to study immersions into $\SS^3$ 
and ${\SS}^4$, in particular the Willmore surfaces. It 
will be done in a separate paper. 

Immersions into the Riemann spaces with non Euclidean 
signature are of 
interest too. For Minkowski space $M^4$: $g_{ik}=diag(1,1,1-1)$. 
The formulae (\ref{eq:3.15}), (\ref{eq:3.16}), (\ref{eq:4.1}), 
(\ref{eq:4.2}) define a surface in $M^4$ with the line element
\beq
ds^2=\left[ 2g_{z\zba}+2 Re \left( g_{zz} \right) \right] dx^2
-4 Im \left(g_{zz} \right) dx dy
+\left[ 2g_{z\zba}-2 Re \left( g_{zz} \right) \right] dy^2
\label{eq:4.8}
\ee
where $z=x+iy$ and
\begin{eqnarray}
g_{zz}=\frac{1}{2} \left( \psiba_1 \varphi_2 -\psiba_2 \varphi_1 
\right)^2=\overline{g_{\zba \; \zba}}
\;\; , \;\;
g_{z \zba}=\frac{1}{2} \left| \psiba_1 \varphi_2 +\psiba_2 
\varphi_1
\right|^2 \;\; .
\label{eq:4.9}
\end{eqnarray}
For a space with the metric $g_{ik}=diag(1,1,-1,-1)$, in 
addition to the immersion of the 
type (\ref{eq:4.8}), (\ref{eq:4.9}) there 
is a Weierstrass type representation for conformal immersion.
\begin{theorem}
 The Weierstrass type formulae 
\begin{eqnarray}
X^1 & = & \frac{i}{2} \int_{\Gamma} { \left[ \left( 
\psiba_1 \psiba_2-\varphi_1 \varphi_2 \right) dz'-
\left( \psi_1 \psi_2 -
\vphiba_1 \vphiba_2 \right) d\zba' \right] }  
\; , \nonumber \\
X^2 & = & \frac{1}{2} \int_{\Gamma} { \left[ \left( 
\psiba_1 \psiba_2+\varphi_1 \varphi_2 \right) dz' +
\left( \psi_1 \psi_2 
+ \vphiba_1 \vphiba_2 \right) d\zba' \right] }  
\; ,\nonumber \\
X^3 & = & -\frac{1}{2} \int_{\Gamma} { \left[ \left( 
\psiba_1 \varphi_2+\psiba_2 \varphi_1 \right) dz'+
\left( \psi_1 \vphiba_2 
+\psi_2 \vphiba_1  \right) d\zba' \right] }  
\; , \nonumber \\
X^4 & = & \frac{i}{2} \int_{\Gamma} { \left[ \left( 
\psiba_1 \varphi_2 - \psiba_2 \varphi_1 \right) dz'-
\left( \psi_1 \vphiba_2 
-\psi_2 \vphiba_1  \right) d\zba' \right] }  
\label{eq:4.10}
\end{eqnarray}
where
\begin{eqnarray}
\begin{array}{l}
\psi_{\alpha z}=p \varphi_\alpha \; , \\
\varphi_{\alpha \zba}=p \psi_\alpha
\end{array}
& \;\; , \;\; \alpha=1,2
\label{eq:4.11}
\end{eqnarray}
and $p$ is a real-valued function, define a conformal 
immersion of surface into the four-dimensional 
space with the metric $g_{ik}=diag(1,1,-1,-1)$. 
The induced metric is of the form
\beq
ds^2=\left( |\psi_1|^2-|\varphi_1|^2 \right) 
\left( |\psi_2|^2-|\varphi_2|^2 \right) dz d\zba= v dz d\zba
\label{eq:4.12}
\ee
and the Gaussian and mean curvatures are
\beq
K=-\frac{2}{v} \left( \log v \right)_{z \zba} \;\; , \;\;
\vec{H}^2=-\frac{4p^2}{v} \;\; .
\label{eq:4.13}
\ee
The total squared mean curvature is
\beq
W=\int {\vec{H}^2 \left[ dS \right] }=-4\int {p^2 dx dy} \;\; .
\label{eq:4.14}
\ee
\end{theorem}
The proof is similar to that the theorem $3.1$. In this case
\begin{eqnarray}
\left( \psi_1 \psi_2 \right)_{z}=
\left( \varphi_1 \varphi_2 \right)_{\zba}
\;\; , \;\;
\left( \psi_1 \vphiba_2 \right)_{z}=
\left( \psiba_2 \varphi_1 \right)_{\zba}
\end{eqnarray}
and
\begin{eqnarray}
d\left( X^1+iX^2 \right) & = & i\psiba_1 \psiba_2 dz +  
i \vphiba_1 \vphiba_2 d\zba \; , \nonumber \\
d\left( X^1-iX^2 \right) & = & -i \varphi_1 \varphi_2 dz  
- i\psi_1 \psi_2 d\zba \; , \nonumber \\
d\left( X^4+iX^3 \right) & = & -i\psiba_2 \varphi_1 dz -
i \psi_1 \vphiba_2 d\zba \; , \nonumber \\
d\left( X^4-iX^3 \right) & = & i\psiba_1 \varphi_2 dz +
i \psi_2 \vphiba_1 d\zba \;\; .
\label{eq:4.17}
\end{eqnarray}
The Weierstrass representation (\ref{eq:4.10})-(\ref{eq:4.14})
could be useful also for the study of $N=2$ superstring
\bib{r30}.

\section{Surfaces in multidimensional spaces and on Lie groups}
\setcounter{equation}{0}

Any solution $(\psi,\varphi)$ of the system (\ref{eq:2.5}) 
gives rise 
via (\ref{eq:2.6}) to the three coordinates $X^1$, $X^2$, $X^3$. 
Given a pair of solutions $(\psi_1,\varphi_1)$ 
and $(\psi_2,\varphi_2)$ one has two possibilities. The 
first is to generate four coordinates via the formula 
(\ref{eq:3.15}), the second is to get six coordinates: 
$X^1$, $X^2$, $X^3$ via (\ref{eq:2.6}) with the solutions 
$(\psi_1,\varphi_1)$ and $X^4$, $X^5$, 
$X^6$ via (\ref{eq:2.6}) using the solutions 
$(\psi_2,\varphi_2)$. In the latter case, one has 
the conformal immersion of a surface into ${\RR}^6$ 
with the induced metric 
\beq
ds^2=\left( u_1^2+u_2^2 \right) dz d{\zba} \;\; . 
\label{eq:5.1}
\ee
Further, introducing four coordinates $X^7$, $X^8$, 
$X^9$, $X^{10}$ via 
(\ref{eq:3.15}) one can get an immersion into 
${\RR}^{10}$. The induced metric in this case is 
\beq
ds^2=\left( u_1^2+u_1 u_2 +u_2^2 \right) dz d{\zba} \;\; .
\label{eq:5.2}
\ee
In such a manner apparently one can get an immersion of a 
surface into the Euclidean space 
${\RR}^N$ with $N=3n+4m$ where $n$ and $m$ are arbitrary 
integers. The corresponding induced metric is of the form
\beq
ds^2=\left( \sum_{\alpha=1}^{n}{u_\alpha^2}+
\sum_{\alpha\neq\beta}^{m} u_\alpha u_\beta \right) dz d{\zba} 
\label{eq:5.3}
\ee
where $m$ is equal to the number of pairs $\alpha$, $\beta$ 
($\alpha \neq \beta$).

The Weierstrass representations for immersion of 
surfaces into the complex 
spaces are defined 
analogously. Let $(\psi_\alpha,\varphi_\alpha)$
and $(\psi_\beta,\varphi_\beta)$ be any two solutions 
of the system 
(\ref{eq:2.5}), {\em i.e.}
\begin{eqnarray}
\begin{array}{l}
\psi_{\alpha z}=p \varphi_\alpha \; ,\\
\varphi_{\alpha \zba}=-p \psi_\alpha
\end{array}
& , & 
\begin{array}{l}
\psi_{\beta z}=p \varphi_\beta \; , \\
\varphi_{\beta \zba}=-p \psi_\beta
\end{array}
\label{eq:5.4}
\end{eqnarray}
where $p$ is a complex valued function. The system 
(\ref{eq:5.4}) implies
\beq
\left( \psi_\alpha \psi_\beta \right)_z=
-\left( \varphi_\alpha \varphi_\beta \right)_{\zba} \;\; .
\label{eq:5.5}
\ee
Consequently, the integral
\beq
X^\gamma=\sum_{\alpha,\beta} {A^\gamma_{\alpha \beta} 
\int_{\Gamma} {\left(\psi_\alpha \psi_\beta d\zba -
\varphi_\alpha \varphi_\beta dz \right)} }
\label{eq:5.6}
\ee
where $A^\gamma_{\alpha\beta}$ are arbitrary constants does 
not depend on the contour of integration $\Gamma$. 
Treating $N$ functions $X^\gamma(z,\zba)$ ($\gamma=1,\dots,N$) 
as the coordinates in ${\CC}^N$, one gets an immersion of 
a surface into ${\CC}^N$. The induced metric is 
given by
\begin{eqnarray}
g_{zz} & = & \left( \sum_{\gamma=1}^{N} \sum_{\alpha \beta} 
{A^\gamma_{\alpha\beta} \psi_\alpha \psi_\beta }\right)^2 \; , 
\nonumber \\
g_{z \zba} & = & -\left( \sum_{\gamma=1}^{N} \sum_{\alpha \beta} 
{A^\gamma_{\alpha\beta} \psi_\alpha \psi_\beta }\right)
\left( \sum_{\gamma=1}^{N} \sum_{\alpha \beta}
{A^\gamma_{\alpha\beta} \varphi_\alpha \varphi_\beta} \right)
\;, \nonumber \\
g_{\zba \; \zba} & = & 
\left( \sum_{\gamma=1}^{N} \sum_{\alpha \beta} 
{A^\gamma_{\alpha\beta} 
\varphi_\alpha \varphi_\beta} \right) ^2
\;\; .
\label{eq:5.7}
\end{eqnarray} 
Choosing $A^\gamma_{\alpha\beta}$ properly, one 
can get a conformal immersion. 

Finally, let us consider a set of solutions $\psi_\alpha^{(i)}$, 
$\varphi^{(i)}_\alpha$ which solve the systems
\begin{eqnarray}
\begin{array}{l}
\psi_{\alpha z}^{(i)}=p^{(i)} \varphi_\alpha^{(i)} \\
\varphi_{\alpha \zba}^{(i)}=-p^{(i)} \psi_\alpha^{(i)}
\end{array}
& , & 
\begin{array}{l}
i=1,\ldots,N \\
\alpha=1,\ldots,M
\end{array}
\label{eq:5.8}
\end{eqnarray}
with different potentials $p^{(i)}$. Since
\beq
\left[ \psi_\alpha^{(i)} \psi_\beta^{(i)} \right]_z=
-\left[ \varphi_\alpha^{(i)} \varphi_\beta^{(i)} \right]_{\zba}
\label{eq:5.9}
\ee
one can define a set of functions $X^{\alpha \beta}$ via
\beq
X^{\alpha \beta}=\sum_{i=1}^{N} {B_i
\int_{\Gamma}^z {\left[ \psi_\alpha^{(i)} \psi_{\beta}^{(i)} 
d\zba -
\varphi_\alpha^{(i)} \varphi_\beta^{(i)} dz \right]} }
\label{eq:5.10}
\ee
where $B_i$ are arbitrary constants. The formula (\ref{eq:5.10}) 
defines, in fact, a matrix $\bf X$ such that 
${\bf X}_{\alpha\beta}= X^{\alpha \beta}$
($\alpha,\beta=1,\ldots ,M$). 
At $N \geq M$ the matrix {\bf $X$} is a 
generic element of the group $GL(M,C)$. So the formula 
(\ref{eq:5.10}) with $N=M$ 
defines a surface on the group $GL(M,C)$ in a meaning given in 
\bib{r31}. Using the formulae from \bib{r31}, one can 
calculate all characteristics of such surfaces.

\section{Integrable deformations}
\setcounter{equation}{0}

An important advantage of the generalized Weierstrass 
representation is 
that it provides a way to construct integrable deformations 
of immersed 
surfaces. The idea is the following 
\bib{r10}-\bib{r11}: let the functions $p$, $\psi$ 
and $\varphi$ in (\ref{eq:2.5}) depend on the parameter $t$. 
Then one 
considers those deformations of $\psi$ and $\varphi$ that 
there are differential operators $A$, $B$, $C$, $D$ such that 
\beq
\begin{array}{l}
\psi_{t}=A\psi+B\varphi \; \; , \\
\varphi_{t}=C\psi+D\varphi \;\; .
\end{array}
\label{eq:6.1}
\ee
Given $A$, $B$, $C$, $D$ the compatibility condition of 
(\ref{eq:6.1}) with (\ref{eq:2.5}) is 
equivalent to the nonlinear partial differential 
equation for $p$. Varying operators 
$A$, $B$, $C$, $D$, one gets an infinite 
hierarchy of integrable equations for $p$. It is the 
so-called modified Veselov-Novikov (mVN) hierarchy \bib{r11}. 
Choosing $A$, $B$, $C$, $D$, as the first 
order operators, one gets the linear equation for $p$. The 
first nontrivial nonlinear example is 
given by the modified Veselov-Novikov equation:
\beq
\begin{array}{lll}
p_{t} & = & p_{zzz}+3p_z \omega 
+\frac{3}{2} p\omega_z +c.c \; , \\
\omega_{\zba} & = & \left( p^2 \right)_z \;\; .
\end{array}
\label{eq:6.2}
\ee
The corresponding operators are
\begin{eqnarray}
A & = & \paz^3 + \pazba^3 + 3 \overline{\omega} \pazba +
\frac{3}{2} \overline{\omega}_{\zba} \;,  \\
B & = & -3p_z \paz + 3p \omega \; , \\
C & = & 3p_{\zba} \pazba - 3p \overline{\omega} \; , \\
D & = & \paz^3 + \pazba^3 + 3 \omega \paz + \frac{3}{2} 
\omega_z \;\; .
\label{eq:6.3}
\end{eqnarray}
The deformation of $\psi$, $\varphi$ via (\ref{eq:6.1}) 
generates the corresponding deformations of the 
coordinates $X^i(z,\zba,t)$. In the case
of immersion into ${\RR}^3$ it was shown in \bib{r14}-\bib{r15}
 that the Willmore functional 
$W$ is invariant under the deformations generated 
by the mVN equation (\ref{eq:6.2}) as well as by the 
whole mVN hierarchy. Formulae for 
deformation of coordinates, elements of metrics and other 
geometric quantities have been obtained in \bib{r19}.

To get integrable deformations of surfaces immersed in 
${\RR}^4$ we assume that both solutions $(\psi_1, \varphi_1)$ 
and $(\psi_2, \varphi_2)$ of the system 
(\ref{eq:3.16}) evolve in $t$ according to equation 
(\ref{eq:6.1}) with the same $A$, $B$, $C$, $D$. 
Correspondingly the coordinates $X^i$ ($i=1,2,3,4$) 
of the surface given by (\ref{eq:3.15}) are deformed 
in $t$. These deformations of a surface are integrable 
one similar to the case ${\RR}^3$ \bib{r11}. 
From (\ref{eq:3.20}) and equality $\int {\left(p^2 \right)_t dz 
d\zba}=0$ it immediately follows 
\begin{theorem}
 The value of the Willmore 
functional $W$ for surface immersed into 
${\RR}^4$ is preserved by the mVN deformations 
(by all hierarchy).
\end{theorem}

There is an infinite set of functionals over $p$ preserved 
by the evolution 
(\ref{eq:6.2}). So there is an infinite family of geometric 
functionals over surface in ${\RR}^4$ which are preserved by 
the mVN deformations. The Willmore functional $W$ for surfaces 
in $R^4$
is invariant under conformal transformations in this space (see 
{\em e.g.} \bib{r22}). One could conjecture that similar to the 
${\RR}^3$ case \bib{r18} all these higher preserved
functionals are invariant 
under the conformal transformations in ${\RR}^4$ too.

The formulae (\ref{eq:6.1}), (\ref{eq:6.2}) define also 
integrable deformations of surfaces in Riemann spaces 
considered above. 
In particular we have 
\begin{theorem} 
The value of Willmore functional for surfaces 
immersed into conformally-Euclidean spaces is 
preserved by the mVN deformations.
\end{theorem} 
Using the results of the papers \bib{r17}, \bib{r20}, \bib{r21}, 
one can define the global deformations of 
surfaces immersed into ${\RR}^4$, ${\SS}^4$ and other spaces. 
This problem will be 
considered elsewhere. 

 \hfill
 \hfill


   \begin{centerline} 
   {\bf REFERENCES}
   \end{centerline}

 \begin{enumerate}

 \item \label{r1}
	G. Darboux, {\it Lecons sur la th\'{e}orie des surfaces 
	et les applications geometriques du calcul 
	infinitesimal}, {\em t. 1-4}, Gauthier-Villars, 
	Paris, 1877-1896.	
  \item \label{r2}
	L. Bianchi, {\it Lezioni di Geometria Differenziale}, 
	{\it 2nd ed.}, Spoerri, Pisa, 1902.
	
  \item \label{r3}
	L.P. Eisenhart, {\it A treatise on the differential 
	geometry of Curves and Surfaces}, Dover, New York, 1909.

  \item \label{r4}
	S.T. Yau (Ed.), {\it Seminar on Differential Geometry}, 
	Princeton Univ., Princeton, 1992.
  \item \label{r5}
	S.S. Chern, in: {\it Differential Geometry and 
	Complex Analysis} (Eds., I. Chavel, H.M. Farkas), 
	Springer, Berlin, 1985.
  \item \label{r6}
	D. Nelson, T. Piran and S. Weinberg (Eds.),
	{\it Statistical mechanics of membranes 
	and surfaces}, World Scientific, Singapore, 1989.
  \item \label{r7}
	D.J. Gross, T. Piran and S. Weinberg (Eds.), {\it Two 
	dimensional quantum gravity and random surfaces}, 
	World Scientific, Singapore, 1992. 
	
  \item \label{r8}
	F. David, P. Ginsparg and Y. Zinn-Justin (Eds.),
	{\it Fluctuating geometries in Statistical Mechanics 
	and Field Theory}, Elsevier Science, Amsterdam, 1996.
  \item \label{r9} 
	K. Kenmotsu, {\it Weierstrass formula for surfaces of 
	prescribed mean curvature}, 
	Math. Ann., {\bf 245}, 89-99, (1979).
  \item \label{r10}
	B.G. Konopelchenko, {\it Multidimensional integrable 
	systems and dynamics of surfaces in space}, preprint of 
	Institute of Mathematics, Taipei, August 1993.
  \item \label{r11}
	B.G. Konopelchenko, {\it Induced surfaces and their 
	integrable dynamics}, Stud. Appl. Math., {\bf 96}, 9-51,
	(1996); preprint Institute of Nuclear Physics, N 93-114, 
	Novosibirsk, (1993).
  \item \label{r12}
	D.A. Hoffman and R. Osserman, {\it The Gauss map of 
	surfaces in ${\RR}^N$}, J. Diff. Geometry, 
	{\bf 18}, 733-754, (1983).
  \item \label{r13}
        D.A. Hoffman and R. Osserman, {\it The Gauss map of 
	surfaces in ${\RR}^3$ and ${\RR}^4$}, 
	Proc. London Math. Soc., (3), {\bf 50}, 27-56, (1985).
  \item \label{r14}
	R. Carroll and B.G. Konopelchenko, {\it Generalized 
	Weierstrass-Enneper inducing, conformal immersion and 
	gravity}, Int. J. Modern Physics {\bf A11}, 
	1183-1216, (1996).
  \item \label{r15}
	B.G. Konopelchenko and I. Taimanov, {\it Generalized 
	Weierstrass formulae, soliton equations and Willmore 
	surfaces}, preprint N. 187, Univ. Bochum, (1995).
  \item \label{r16}
        B.G. Konopelchenko and I. Taimanov, {\it Constant mean 
	curvature surfaces via an integrable dynamical system}, 
	J. Phys. {\bf A}:Math. Gen., {\bf 29}, 1261-1265, (1996).
 \item \label{r17}
	I. Taimanov, {\it Modified Novikov-Veselov equation and 
	differential geometry of surfaces}, Trans. 
	Amer. Math. Soc., Ser. 2, {\bf 179}, 133-159, (1997).
 \item \label{r18}
	P.G. Grinevich and M.V. Schmidt, {\it Conformal 
	invariant functionals of immersion of tori into 
	${\RR}^3$}, Journal of Geometry and Physics (to appear); 
	preprint SFB288 N 291, TU-Berlin, 1997.
 \item \label{r19}
	J. Richter, {\it Conformal maps of a Riemann surface 
	into space of quaternions}, PH.D Thesis, TU-Berlin, 1997.
 \item \label{r20}
        I. Taimanov, {\it Global Weierstrass representation and 
	its spectrum}, Uspechi Mat. Nauk, {\bf 52}, N 6, 
	187-188, (1997).
 \item \label{r21}
        I. Taimanov, {\it The Weierstrass representation of 
	spheres in ${\RR}^3$, the Willmore 
	numbers and soliton spheres}, preprint SFB 288, 
	N 302, TU-Berlin, 1998.
 \item \label{r22}
        T.J. Willmore, {\it Total curvature in Riemannian 
	Geometry}, Ellis Horwood, New York, 1982.
 \item \label{r23}
	K. Nomizu and T. Sasaki, {\it Affine Differential 
	Geometry}, Cambridge Press, 1994.
 \item \label{r24}
	U. Abresch, {\it Spinor representation of CMC surfaces}, 
	Lecture at Luminy, 1989.
 \item \label{r25}
	A.I. Bobenko, {\it Surfaces in terms of $2$ by $2$ 
	matrices. Old and new integrable cases}, in: {\it 
	Harmonic maps and integrable 
	systems}, (A. Fordy and J. Wood, Eds.), Vieweg, 1994.
 \item \label{r26}
	R. Kusner and N. Schmitt, {\it The spinor representation 
	of surfaces in space}, preprint dg-ga/9610005, 1996.
 \item \label{r27}
	T. Friedrich, {\it On the spinor representation of 
	surfaces in euclidean $3-$spaces}, 
	preprint SFB 288, N 295, TU-Berlin, 1997. 
 \item \label{r28}
	A.I. Bobenko, {\it All constant mean curvature tori in 
	${\RR}^3$, ${\SS}^3$, $\HH^3$ 
	in terms of theta-functions}, Math. Ann., {\bf 290},
	 209-245, (1991).
 \item \label{r29}
	L.P. Eisenhart, {\it Riemannian Geometry}, 
	Princeton, 1926.
 \item \label{r30}
	J. Barrett, G.W. Gibbons, M.J. Perry and P. Ruback, 
	{\it Kleinian Geometry and the $N=2$ superstring}, 
	Int. J. Modern Phys. A, {\bf 9}, 1457-1493, (1994).
 \item \label{r31}
	A.S. Fokas and I.M. Gelfand, {\it Surfaces on Lie 
	groups, on Lie algebras and their integrability}, 
	Commun. Math. Phys., {\bf 177}, 203-220, (1996).

\end{enumerate}
\end{document}